\documentclass[12pt]{article}

\usepackage{latexsym}
\usepackage{amsmath,amssymb,amsthm}

\usepackage{amsmath,amsthm,latexsym,amssymb,pgf,makeidx,
dsfont,authblk,boxedminipage,tikz}
\usetikzlibrary{arrows,decorations.pathmorphing,backgrounds,positioning,fit,petri,calc}
\date{}

\begin{document}
	
	\newtheorem{definition}{Definition}
	\newtheorem{theorem}{Theorem}
	\newtheorem{corollary}{Corollary}
	\newtheorem{proposition}{Proposition}
	\newtheorem{rem}{Rem:}
	\newtheorem{lemma}{Lemma}
	\newtheorem{example}{Example}
	\newtheorem{conjecture}{Conjecture}

	\title{Average eccentricity, minimum degree and maximum degree in graphs}
	
\author{P. Dankelmann, F.J. Osaye 
	(University of Johannesburg) }

\maketitle

\begin{abstract} 
Let $G$ be a connected finite graph with vertex set $V(G)$. The eccentricity $e(v)$ of 
a vertex $v$ is the distance from $v$ to a vertex farthest from $v$. The average 
eccentricity of $G$ is defined as $\frac{1}{|V(G)|}\sum_{v \in V(G)}e(v)$. 
We show that the average eccentricity of a connected graph of order $n$, minimum degree 
$\delta$ and maximum degree $\Delta$ does not exceed 
$\frac{9}{4} \frac{n-\Delta-1}{\delta+1} 
	     \big( 1 + \frac{\Delta-\delta}{3n} \big) + 7$, and this bound is sharp
apart from an additive constant. We give improved bounds for triangle-free graphs
and for graphs not containing a $4$-cycles.  
\end{abstract}

\section{Introduction}

\noindent
Let $G$ be a connected graph with vertex set $V(G)$. The {\em eccentricity} $e(v)$ of a vertex is the maximum distance between $v$ and a vertex in $G$. 
The {\em average eccentricity} ${\rm avec}(G)$ of $G$ is the average of the eccentricities 
of the vertices of $G$, i.e., ${\rm avec}(G) = \frac{1}{|V(G)|}\sum_{u\in V(G)} e(u)$. The 
average eccentricity was first introduced by Buckley and Harary \cite{BucHar} under the 
name eccentric mean. Its systematic study was initiated by Dankelmann, Goddard and Swart 
\cite{DanGodSwa2004},  who established an upper bound on the average eccentricity  in terms 
of order and minimum degree, and further showed among other results, that the path 
maximizes the average eccentricity among all connected graphs of given order. 
\begin{proposition}[\cite{DanGodSwa2004}]\label{prop1}
Let $G$ be a connected graph of order $n$. Then 
\begin{equation*}
	{\rm avec(G)} \leq {\rm avec}(G) \leq \frac{1}{n} \left\lfloor \frac{3n^2}{4} -\frac{n}{2} \right\rfloor,  	
\end{equation*}
and this bound is sharp.
\end{proposition}
We first recall some upper bounds on the average eccentricity of a connected graph 
involving the minimum degree. Dankelmann, Goddard and Swart \cite{DanGodSwa2004} 
showed that if $G$ is a graph of order $n$ and minimum degree $\delta \geq 2$, then
\begin{equation}\label{equ-DGS}
{\rm avec}(G) \leq \frac{9n}{4(\delta +1)} + \frac{15}{4},
\end{equation}
and this bound is sharp apart from a small additive constant. By applying a similar technique, it was recently proved in \cite{DanMukOsaRod} that if $G$ is $K_3$-free, then the bound \eqref{equ-DGS} can be improved to 
\begin{equation}\label{eqn-DMOR-K3}
{\rm avec}(G) \leq 3\left\lceil \frac{n}{2\delta} \right\rceil +5,
\end{equation}
and for $C_4$-free graphs to 
\begin{equation}\label{eqn-DMOR-C4}
	{\rm avec}(G) \leq \frac{15}{4} \left\lceil \frac{n}{\delta^2 -2\lfloor \frac{\delta}{2} \rfloor +1} \right\rceil +\frac{11}{2}.
\end{equation}
Moreover, they showed that the bound in \eqref{eqn-DMOR-K3} is sharp apart from a small 
additive constant and that, for $\delta +1$ a prime power, there exists an infinite number 
of $C_4$-free graphs of minimum degree at least $\delta$ with 
\begin{equation*}
	{\rm avec}(G) \geq \frac{15}{4} \frac{n}{\delta^2 +3\delta +2} -\frac{5}{2}.
\end{equation*}
Other recent results on the average eccentricity of graphs can be found, for example, in 
 (\cite{AliDanMorMukSwaVet2018}, \cite{DanMuk2014}, \cite{DanOsa-manu}, \cite{DanMukOsaRod}, 
 \cite{DarAliKlaDas2018}, 
  \cite{DuIli2013}, \cite{DuIli2016}, \cite{Ili2012}, \cite{SmiSzeWan2016}, \cite{TanZho2012}).

In the bounds \eqref{equ-DGS}, \eqref{eqn-DMOR-K3} and \eqref{eqn-DMOR-C4}, the examples that 
show that these bounds are sharp or close to being sharp, all have vertex degrees close to the 
minimum degree. That suggests that these bounds can be improved for graphs containing a vertex 
of large degree. In this paper, we show that this is indeed the case by proving asymptotically 
sharp upper bounds on the average eccentricity of a graph of given order, minimum degree and 
maximum degree. We also give corresponding bounds for triangle-free graphs and $C_4$-free 
graphs. Our bounds improve on the inequalities \eqref{equ-DGS}, \eqref{eqn-DMOR-K3} 
and \eqref{eqn-DMOR-C4}. Moreover, we construct graphs to show that our bounds for connected 
graphs and triangle-free graphs are sharp apart from a small additive constant, and that for 
$C_4$-free graphs our bound is close to being best possible.

\section{Definitions, notation and preliminaries}

\noindent
We use the following notation. Let $G$ be a connected graph of order $n$ with vertex set 
$V(G)$, edge set $E(G)$. The \textit{distance} between two vertices $u$ and $v$, denoted 
by $d_G(u,v)$ is the length of a shortest $(u,v)$-path in $G$.  
The diameter ${\rm diam}(G)$ and radius ${\rm rad}(G)$ are 
the largest and smallest, respectively, of all eccentricities of vertices of $G$. The {\em 
total eccentricity} $EX(G)$ is the sum of the eccentricities of the vertices of $G$.
If $x$ is a vertex of $G$, $e$ an edge of $G$, and $A$ a subset of $V(G)$, then
$d_G(x, A)$ is the minimum of the distances between $x$ and a vertex in $A$, and 
$d_G(e,A)$ is the minimum of the distances between a vertex incident with $e$ and the 
vertices in $A$. 

The neighbourhood of a vertex $v$, denoted by $N_G(v)$, is the set of vertices adjacent to $v$. The closed neighbourhood $N_G[v]$ is defined as $N_G(v) \cup \{v\}$. For 
$k\in \mathbb{N}$, the $kth$ neighbourhood of a given subset $A\subset V$, denoted by 
$N_G^k(A)$, is the set of all vertices $x$ of $G$ with $d_G(x,a) \leq k$ for some $a\in A$. 
If $k=1$, then we simply write $N_G[A]$. If 
no confusion can occur, we drop the argument or subscript $G$. The degree 
${\rm deg}(v)$ of a vertex $v$ is defined  as $|N_G(v)|$. We denote the minimum and 
maximum degree of $G$ by $\delta(G)$ and 
$\Delta(G)$, respectively. For 
$k\in \mathbb{N}$, the $k^{th}$ power $G^k$ of $G$ is the graph with the same vertex set of $G$, in which two distinct vertices $u$ and $v$ are adjacent if $d(u,v)\leq k$. If $A\subset V$, then $G^k[A]$ is the subgraph of $G^k$ induced by $A$. If $H$ is a subgraph of $G$, we write $H\leq G$.

The line graph $L(G)$ of $G$ is the graph whose vertices are the edges of $G$ such that two vertices of $L(G)$ are adjacent if they share a vertex as edges of $G$. Let $M\subseteq E(G)$, then $V(M)$ is the set of vertices incident with at least an edge of $M$.  For $k\in \mathbb{N}$, a $k$-packing of $G$ is a subset $A\subseteq V$ with $d_G(u,v) >k$ for all $u,v \in A$. The sequential sum $G_1 +G_2 +\ldots +G_k$ is the graph obtained from the disjoint union of the graphs $G_1$, $G_2, \ldots, G_k$ by joining every vertex of $G_i$ to every vertex of $G_{i+1}$ for $i=1,2,\ldots, k-1$.

A set of edges $M\subset E(G)$ is a \textit{matching} in $G$ if no two edges of $M$ are 
incident.  A \textit{maximum matching} is a matching of maximum size. A graph $G$ is said to be \textit{complete} if all vertices of $G$ are pairwise adjacent. A path, cycle or complete graph of order $n$ is denoted by $P_n$, $C_n$ or $K_n$. We refer to $K_3$ also as a \textit{triangle}. For a positive integer $k$, $kK_1$ is an edgeless graph of $k$ isolated vertices. If $F$ is a graph, then $G$ is said to be $F$-free if $G$ does not contain $F$
as a subgraph. 

An important tool in the proofs of our main results is the weighted eccentricity of a graph defined in \cite{DanGodSwa2004}.
\begin{definition}[\cite{DanGodSwa2004}]\label{def1}
	Let $G$ be a connected graph and $c:V(G) \rightarrow  \mathbb{R}$ a 
	nonnegative weight function
	on the vertices of $G$. Then the eccentricity of $G$ with respect to $c$ is defined by 
	\[ EX_c(G)= \sum_{x\in V(G)} c(x) e_G(x). \] 
	Let $N=\sum_{x\in V(G)} c(x)$ be the total weight of the vertices in $G$. If $N>0$, then the average eccentricity of $G$ with respect to $c$ is 
	\[ {\rm avec}_c(G) = \frac{EX_c(G)}{N}.   \]
\end{definition}
The following result from \cite{DanGodSwa2004} generalises the bound in Proposition \ref{prop1}.
\begin{proposition}[\cite{DanGodSwa2004}]\label{prop2}
	Let $G$ be a connected graph, $c$ a weight function on the vertices 
	of $G$, and $N = \sum_{v\in V(G)} c(v)$ the total weight of the vertices
	of $G$. If $c(v) \geq 1$ for all $v\in V(G)$, then 
	\[
	{\rm avec}_c(G) \leq {\rm avec}(P_{\lceil N \rceil}).
	\]
\end{proposition}

\section{Main Results}

\noindent
We begin by presenting a result on an upper bound on the average eccentricity in terms of order, minimum degree and maximum degree. The approach used throughout this section is a slight modification of that used in \cite{DanGodSwa2004} and \cite{DanMukOsaRod}.

\subsection{A bound on the average eccentricity of a graph}
\begin{theorem}\label{thm5-1}
	Let $G$ be a connected graph of order $n$, minimum degree $\delta$ and maximum degree $\Delta$. Then,
	\begin{equation*}\label{label1}
	{\rm avec}(G)\leq   \frac{9}{4} \frac{n-\Delta-1}{\delta+1} 
	     \Big( 1 + \frac{\Delta-\delta}{3n} \Big) + 7. 
	\end{equation*}
	This bound is sharp apart from the value of the additive constant.
\end{theorem} 

\noindent
{\bf Proof:}
Let $v_1$ be a vertex of degree $\Delta$. We find a maximal 2-packing $A$ of $G$ using the 
following method: Let $A=\{v_1\}$. If there exists a vertex $a_1$ with $d_G(a_1,A)=3$, add 
$a_1$ to $A$. Add vertices $a_i$ with $d_G( a_i,A)=3$ to $A$ until every vertex not in $A$ is 
within distance two of $A$. Then $A=\{v_1,a_1,a_2,\ldots, a_r\}$ and $|A|=r+1$. Let $T_1$ be 
the forest with vertex set $N[A]$ and whose edge set consists of all edges incident with 
vertices in $A$. 

By our construction of $A$, there exist $r$ edges in $G$, each joining two neighbours of 
distinct vertices of $A$, such that the addition of these edges to $T_1$ yields a subtree 
$T_2 \leq G$. Now every vertex $u\in V(G)-V(T_2)$ is adjacent to some $u'\in V(T_2)$. Let 
$T$ be the spanning tree of $G$ with edge set $E(T) = E(T_2) \cup \{uu'| u\in V(G)-V(T_2)\}$. 
Since ${\rm deg}_T(v_1) = {\rm deg}_G(v_1)$, it follows that $T$ has the same maximum 
degree as $G$.

Since deleting edges does not decrease the average eccentricity, we have 
${\rm avec}(G) \leq {\rm avec}(T)$.  Thus it suffices to prove the result for $T$.  
We think of $T$ as a weighted tree in which each vertex has weight $1$. We now obtain a 
new weight function by moving the weight 
every vertex to a nearest vertex in $A$. More specifically, for every vertex $x$ of $T$ 
let $x_A$ be a vertex in $A$ closest to u in $T$. 
Let $c:V(T)\rightarrow \mathbb{N} \cup \{0\}$ be the weight function defined by 
\[c(u)=|\{x\in V(T)|x_A=u\}|\] for each $x\in V(T)$. Then $c(u)=0$ if $u\not \in A$, $c(u) \geq  \delta +1$ if $u\in A-\{v_1\}$, and $c(u) \geq \Delta +1$ if $u = v_1$. Note that $\sum_{u\in V(T)} c(u) =n$. It follows that 
\[
n = \sum_{u\in V(T)} c(u) = \sum_{u\in A}c(u) 
 \geq \Delta +1 + (|A|-1)(\delta +1)
 =|A|(\delta +1)+\Delta -\delta.
\]
By rearranging, we have 
\begin{equation}\label{eqn5.1*}
|A| \leq \frac{n-\Delta +\delta}{\delta +1}.
\end{equation}
For every vertex $x$ of $T$ we have $d(x,x_A) \leq 2$ and thus $|e_T(x) - e_T(x_A)| \leq 2$. Hence 
\begin{eqnarray}
{\rm avec}(T) & = & \frac{1}{n} \sum_{x \in V(T)} e_T(x) \nonumber \\
& \leq & \frac{1}{n} \sum_{x\in V(T)} (e_T(x_A) + 2) \nonumber \\
& = & \big( \frac{1}{n} \sum_{u \in A} c(u) e_T(u) \big) + 2 \nonumber \\
& = & {\rm avec}_{c}(T) + 2.    \label{eqn-5b} 
\end{eqnarray}
Note that the weight of $c$ is concentrated in $A$. Thus we consider the induced subgraph 
$T^3[A]$, which by the construction of $A$ is connected. 
Clearly for all $a_i,a_j \in A$, $d(a_i,a_j)\leq 3d_{T^3[A]}(a_i,a_j)$. Since every vertex 
of $T$ is within distance two of some vertex in $A$, we have for each $a_i\in A$, 
$e_{T}(a_i) \leq 3e_{T^3[A]}(a_i) +2$. Thus, 
\begin{equation}\label{eqn5a}
{\rm avec_c}(T) \leq 3{\rm avec_c}(T^3[A]) +2.
\end{equation}
We now modify the weight function $c$ to obtain a new weight function $c'$
which satisfies $c'(a) \geq 1$ for all $a\in A$. 
Define $c'(u) = \frac{c(u)}{\delta+1}$ for $u\in A-\{v_1 \}$, and  
$c'(v_1) = \frac{c(v_1)-\Delta+\delta}{\delta+1}$. 
Since $c(u) \geq {\rm deg}_G(u)+1$ for all $u\in A$, and since ${\rm deg}_G(v_1)=\Delta$ 
while ${\rm deg}_G(a_i) \geq \delta$ for $a_i \in A - \{v_1\}$, we have
$c'(u) \geq 1$ for all $u\in A$. 
 
Since  $\sum_{u\in A}c(u) = n$, we have for  the total weight $N$ of $c'$, 
\[N= \sum_{u\in A}c'(u)=\frac{(\sum_{u\in A}c(u))-\Delta +\delta}{\delta +1}
=  \frac{n-\Delta +\delta }{\delta + 1}.  \] 
We now express ${\rm avec_c}(T^3[A])$ in terms of ${\rm avec}_{c'}(T^3[A])$. We have 
\begin{eqnarray*}
	{\rm avec}_{c'}(T^3[A]) &=&\frac{EX_{c'}(T^3[A])}{N}\\
	&=& \frac{\frac{1}{\delta +1} [\sum_{u\in A-\{v_1\}}c(u) e_{T^3[A]}(u)+ (c(v_1)-\Delta +\delta)e_{T^3[A]}(v_1) ]}{\frac{1}{\delta +1} [n-\Delta +\delta]}   \\
	&=& \frac{EX_{c}(T^3[A])- e_{T^3[A]}(v_1)(\Delta -\delta)}{n-\Delta+\delta} \\
	&=& \frac{n}{n-\Delta+\delta} {\rm avec}_{c}(T^3[A]) 
	- \frac{\Delta-\delta}{n-\Delta+\delta}  e_{T^3[A]}(v_1), 
\end{eqnarray*}
and thus, by rearranging,  
\begin{equation} \label{eq:avec(Tc)-in-terms-of-avec(Tc')}
{\rm avec_c}(T^3[A]) = \frac{n-\Delta+\delta}{n} {\rm avec_{c'}}(T^3[A]) 
+\frac{\Delta -\delta}{n} e_{T^3[A]}(v_1).
\end{equation}
We bound the two terms of the right hand side of \eqref{eq:avec(Tc)-in-terms-of-avec(Tc')} 
separately. Since $T^3[A]$ has order $|A|$, and since $|A| \leq \frac{n-\Delta+\delta}{\delta+1}$   
by \eqref{eqn5.1*}, we have
\begin{equation}  \label{eq:bound-on-e(v_1)} 
e_{T^3[A]}(v_1) \leq {\rm diam}(T^3[A]) \leq |A|-1 \leq \frac{n-\Delta-1}{\delta+1}. 
\end{equation}
To bound ${\rm avec_{c'}}(T^3[A])$ note that $c'(u) \geq 1$ for all $u\in A$. 
Hence we have by Proposition \ref{prop2}, 
\begin{equation}\label{eqn-5a}
{\rm avec_{c'}}(T^3[A]) 
\leq {\rm avec}(P_{\lceil N\rceil}) 
\leq \frac{3}{4}\lceil N\rceil-\frac{1}{2} 
< \frac{3}{4} \frac{n - \Delta-1}{\delta+1} + 1,
\end{equation}
with the last inequality holding since 
$\lceil N \rceil = \lceil \frac{n-\Delta+\delta}{\delta+1} \rceil < \frac{n-\Delta-1}{\delta+1} + 2$.  
Substituting  \eqref{eq:bound-on-e(v_1)} and \eqref{eqn-5a} into 
\eqref{eq:avec(Tc)-in-terms-of-avec(Tc')} yields, after simplification, 
\begin{eqnarray}
{\rm avec_{c}}(T^3[A])  & < & 
\frac{n-\Delta+\delta}{n} \Big( \frac{3}{4} \frac{n - \Delta-1}{\delta+1} + 1 \Big) 
+\frac{\Delta -\delta}{n} \Big( \frac{n-\Delta-1}{\delta+1} \Big) \nonumber \\ 
& = & \frac{3}{4} \frac{n-\Delta-1}{\delta+1} 
  \big[ 1 + \frac{\Delta-\delta}{3n} \big] + 1 -  \frac{\Delta-\delta}{n} \nonumber \\
& \leq & \frac{3}{4} \frac{n-\Delta-1}{\delta+1} 
  \big[ 1 + \frac{\Delta-\delta}{3n} \big] + 1.  
\label{eqn5-3a}
\end{eqnarray}
Combining inequalities \eqref{eqn-5b}, \eqref{eqn5a}, and \eqref{eqn5-3a}, we obtain  
\begin{eqnarray*}
{\rm avec}(T) 
    &\leq & {\rm avec_{c}}(T) +2 \\
	&\leq &  3 {\rm avec_{c}}(T^3[A]) +4 \\
	& < &  \frac{9}{4} \frac{n-\Delta-1}{\delta+1} 
	     \Big( 1 + \frac{\Delta-\delta}{3n} \Big) + 7.         
\end{eqnarray*}
The bound on ${\rm avec}(G)$ now follows since ${\rm avec}(G) \leq {\rm avec}(T)$. 

To see that the bound in Theorem \ref{thm5-1} is sharp apart from an additive constant, let $\delta$, $\Delta$ and $k$ be positive integers with $\Delta \geq \delta +1$ and consider the graph  $G_{\delta,\Delta,k}$ obtained as follows.   

\begin{figure} 
\begin{center}
	\begin{tikzpicture}
	[scale=0.7,inner sep=0.8mm, 
	vertex/.style={circle,thick,draw}, 
	thickedge/.style={line width=3pt}] 
	\begin{scope}[>=]   
	\node[vertex] (a0) at (-2.5,0) [fill=black] {};
	\node[vertex] (a1) at (-1,0) [fill=black] {};
	\node[vertex] (a2) at (0.5,0) [fill=black] {};    
	\node[vertex] (a3) at (2,0) [fill=black] {};    
	\node[vertex] (a4) at (3.5,0) [fill=black] {};        
	\node[vertex] (a5) at (5,0) [fill=black] {};
	\node[vertex] (a6) at (6.5,0) [fill=black] {};
	\node[vertex] (a7) at (8,0) [fill=black] {};
	\node[vertex] (a8) at (9.5,0) [fill=black] {}; 
	\node[vertex] (a9) at (11,0) [fill=black] {};

	\node[vertex] (b1) at (-1.75,1) [fill=black] {}; 
	\node[vertex] (c1) at (-1.75,-1) [fill=black] {};
	\node[vertex] (b2) at (1.25,1) [fill=black] {}; 
	\node[vertex] (c2) at (1.25,-1) [fill=black] {}; 
	\node[vertex] (b3) at (4.25,1) [fill=black] {}; 
	\node[vertex] (c3) at (4.25,-1) [fill=black] {};
	\node[vertex] (b4) at (7.25,1) [fill=black] {}; 
	\node[vertex] (c4) at (7.25,-1) [fill=black] {};
	\node[vertex] (b5) at (10.25,1) [fill=black] {}; 
	\node[vertex] (c5) at (10.25,-1) [fill=black] {};
	
	\node[vertex] (d1) at (12.5,0) [fill=black] {};
	\draw[black, very thick] (d1)--(a9);
	\node[vertex] (d2) at (14,-1) [fill=black] {};
	\node[vertex] (d3) at (14,1) [fill=black] {};
	\node[vertex] (d4) at (15.5,-2) [fill=black] {};
	\node[vertex] (d5) at (15.5,2) [fill=black] {};
	\node[vertex] (d6) at (17,-1) [fill=black] {};
	\node[vertex] (d7) at (17,1) [fill=black] {};
	\node[vertex] (d8) at (18.5,0) [fill=black] {};
	
	\draw[black, very thick] (d1)--(d2);
	\draw[black, very thick] (d1)--(d3);
	\draw[black, very thick] (d1)--(d4);
	\draw[black, very thick] (d1)--(d5);
	\draw[black, very thick] (d1)--(d6);
	\draw[black, very thick] (d1)--(d7);
	\draw[black, very thick] (d1)--(d8);
	
	\draw[black, very thick] (d2)--(d3);
	\draw[black, very thick] (d2)--(d4);
	\draw[black, very thick] (d2)--(d5);
	\draw[black, very thick] (d2)--(d6);
	\draw[black, very thick] (d2)--(d7);
	\draw[black, very thick] (d2)--(d8);
	
	\draw[black, very thick] (d3)--(d4);
	\draw[black, very thick] (d3)--(d5);
	\draw[black, very thick] (d3)--(d6);
	\draw[black, very thick] (d3)--(d7);
	\draw[black, very thick] (d3)--(d8);
	
	\draw[black, very thick] (d4)--(d5);
	\draw[black, very thick] (d4)--(d6);
	\draw[black, very thick] (d4)--(d7);
	\draw[black, very thick] (d4)--(d8);
	
	\draw[black, very thick] (d5)--(d6);
	\draw[black, very thick] (d5)--(d7);
	\draw[black, very thick] (d5)--(d8);
	
	\draw[black, very thick] (d6)--(d7);
	\draw[black, very thick] (d6)--(d8);
	
	\draw[black, very thick] (d7)--(d8);
	
	\draw[black, very thick] (a0)--(a1);
	\draw[black, very thick] (a1)--(a2); 
	\draw[black, very thick] (a3)--(a4);  
	\draw[black, very thick] (a5)--(a6);
	\draw[black, very thick] (a7)--(a8);
	
	\draw[black, very thick] (a0)--(b1);
	\draw[black, very thick] (a0)--(c1);
	\draw[black, very thick] (b1)--(a1);
	\draw[black, very thick] (c1)--(a1);
	\draw[black, very thick] (b1)--(c1);
	
	\draw[black, very thick] (a2)--(b2);
	\draw[black, very thick] (a2)--(c2);
	\draw[black, very thick] (b2)--(a3);
	\draw[black, very thick] (c2)--(a3);
	\draw[black, very thick] (b2)--(c2);
	
	\draw[black, very thick] (a4)--(b3);
	\draw[black, very thick] (a4)--(c3);
	\draw[black, very thick] (b3)--(a5);
	\draw[black, very thick] (c3)--(a5);
	\draw[black, very thick] (b3)--(c3);
	
	\draw[black, very thick] (a6)--(b4);
	\draw[black, very thick] (a6)--(c4);
	\draw[black, very thick] (b4)--(a7);
	\draw[black, very thick] (c4)--(a7);
	\draw[black, very thick] (b4)--(c4);
	
	\draw[black, very thick] (a8)--(b5);
	\draw[black, very thick] (a8)--(c5);
	\draw[black, very thick] (b5)--(a9);
	\draw[black, very thick] (c5)--(a9);
	\draw[black, very thick] (b5)--(c5);
	
	\end{scope}         
	\end{tikzpicture} 
\end{center}
\caption{$G_{\delta,\Delta, k}$ with $\delta =3$, $\Delta =8$ and $k=6$. }
  \label{fig:sharpness-example-for-general-graphs} 
\end{figure}
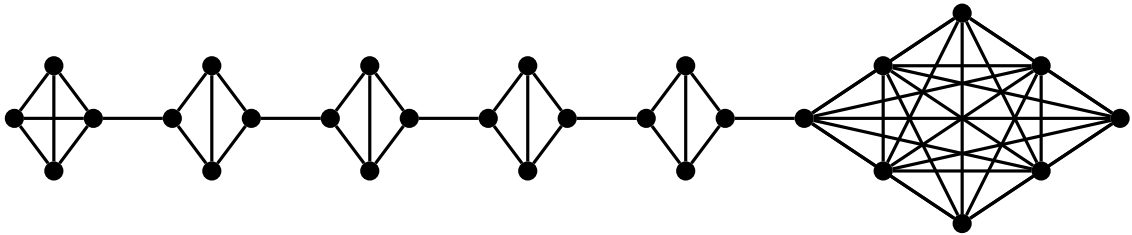

Let $G_1,G_2,\ldots , G_{k-1}$ be disjoint copies of the complete graph $K_{\delta +1}$ and let $G_k$ be the complete graph $K_{\Delta}$. For $i=2,3,\ldots,k-1$ let 
$u_iv_i$ be an edge of $G_i$. Let $G_{\delta,\Delta,k}$ be the graph with vertex set 
$V(G_{\delta,\Delta,k})=V(G_1) \cup V(G_2)\cup \ldots \cup V(G_k)$, 
and edge set 
$E(G_{\delta,\Delta, k})= E(G_1) \cup E(G_2) \cup \cdots \cup E(G_k) 
  - \{u_2v_2, u_3v_3,\ldots,u{k-1}v_{k-1}\} \cup \{v_1u_2, v_2u_3,\ldots, v_{k-1}u_k\}$. 
The graph $G_{3,8,6}$ is shown in  Figure \ref{fig:sharpness-example-for-general-graphs}. 
  
It is easy to see that $G_{\delta,\Delta,k}$ has minimum degree $\delta$, 
maximum degree $\Delta$, order $n=\Delta+(\delta +1)(k-1)$,  
${\rm diam}(G_{\delta,\Delta, k}) = 3(k-1)$ and 
${\rm rad}(G_{\delta,\Delta, k}) = \lceil \frac{3}{2}(k-1)\rceil$. 
A straightforward calculation, whose details we omit, shows that
\[ EX(G_{\delta,\Delta,k}) =
  3(k-1)(\Delta+\delta+1) - 2 
   + 6(\delta+1)[\frac{3}{8}k^2 - \frac{5}{4}k + 1]. \]
   
Now $k-1=\frac{n-\Delta}{\delta+1}$. Hence  
$3(k-1)(\Delta+\delta+1) = 3 \frac{n-\delta}{\delta +1}(\Delta+ \delta+1)
       =  3\frac{(n-\Delta)\Delta}{\delta+1} + O(n)$ 
and
$6(\delta+1)[\frac{3}{8}k^2 - \frac{5}{4}k + 1] = 6(\delta+1)[\frac{3}{8}(k-1)^2 +O(k)] 
     = \frac{9}{4} \frac{(n-\Delta)^2}{\delta+1} +O(n)$,
and thus 
\begin{eqnarray*} 
{\rm avec}(G_{\delta,\Delta,k})    
 &= &\frac{1}{n} EX(G_{\delta,\Delta,k})  \\
 & = &  \frac{1}{n} 
    \Big( 3 \frac{(n-\Delta)\Delta}{\delta+1} 
    + \frac{9}{4} \frac{(n-\Delta)^2}{\delta+1} +O(n) \Big)   \\
 & = & \frac{9}{4} \frac{n-\Delta}{\delta+1} (1 + \frac{\Delta}{3n}) + O(1). 
\end{eqnarray*} 
Since the upper bound on ${\rm avec}(G)$ is also 
$\frac{9}{4} \frac{n-\Delta}{\delta+1} (1 + \frac{\Delta}{3n}) + O(1)$,
we conclude that the bound is sharp apart from an additive constant. 
\qed \\

We note that for fixed $\delta$ and large $n$, the value $avec(G_{\delta,\Delta,k})$ and 
the bound in Theorem \ref{thm5-1} differ only by $O(1)$, irrespective of whether $\Delta$ is 
constant or grows with $n$. In addition, since $\Delta \geq \delta$ and since the bound 
in Theorem \ref{thm5-1} is decreasing in $\Delta$, replacing 
$\Delta$ by $\delta$ in our bound in Theorem \ref{thm5-1} yields a bound on the 
average eccentricity in terms of order and minimum degree that differs from  
\eqref{equ-DGS} only by 
a small additive constant. Hence Theorem \ref{thm5-1} is, in some sense, a generalisation 
of \eqref{equ-DGS}.

\subsection{A bound on the average eccentricity of a triangle-free graph}

We now show that the bound in Theorem \ref{thm5-1} can be improved for triangle-free graphs. 

\begin{theorem}\label{thm5-2}
	Let $G$ be a connected triangle-free graph with $n$ vertices, minimum degree $\delta$ and maximum degree $\Delta$. Then, 
	\begin{equation*}\label{eqn5Dd}
	{\rm avec} (G) \leq \frac{3}{2} \frac{n-\Delta}{\delta} 
    \Big[ 1 + \frac{\Delta-\delta}{3n} \Big] 
    + \frac{19}{2}.
	\end{equation*}
	This bound is sharp apart from the value of the additive constant.
\end{theorem}

\noindent
\textbf{Proof.}
Let $v_1$ be a vertex of degree $\Delta$ in $G$ and let $e_1$ be an edge incident with $v_1$. 
We obtain a maximal matching $M$ of $G$ as follows.
Let $M=\{ e_1\}$. Let $V(M)$ be the set of vertices incident with an edge of $M$. Recall that 
$d_G(e_1,V(M))$ is the minimum of the distances between a vertex incident with edge $e_1$ and 
a vertex in $V(M)$.
If there exists an edge $e_2$ in $G$ with $d_G(e_2,V(M))=3$, add $e_2$ to $M$. Add edges $e_i$ 
with $d_G(e_i,V(M))=3$ to $M$ until each of the edges not in $M$ is within distance two of 
$M$. 

Let $T_1\leq G$ be a subforest of $G$ with vertex set $N[V(M)]$ and edge set consisting of all 
edges incident with a vertex in $V(M)$. By our construction of $M$, $G$ contains $|M|-1$ 
edges, each joining two distinct components of $T_1$, whose addition to $T_1$ yields a 
subtree $T_2 \leq G$ so that $T_2$ contains $T_1$ and $|V(T_1)| = |V(T_2)|$. 
 
Now every vertex $x$ not in $T_2$ is within distance at most three from some vertex $x'$ in 
$T_2$. Let $T$ be a spanning tree of $G$ which contains $T_2$ and which is distance preserving 
from $V(M)$, i.e. $d_T(x,V(M))=d_G(x,V(M))$, for every vertex $x\in V(G)$. Since 
${\rm deg}_T(v_1) = {\rm deg}_G(v_1)$, tree $T$ has the same maximum degree as $G$. 
Furthermore, since 
${\rm avec}(G)\leq {\rm avec}(T)$, it suffices to prove the bound for $T$. 
 
For every vertex $u \in V(T)$, let $u_M$ be a vertex in $V(M)$ closest to $u$ in $T$. We 
can view $T$ as a weighted tree in which each vertex has weight $1$. We now move the weight 
of $u$ to $u_M$. That is, we define a weight function 
$c:V(T) \rightarrow \mathbb{N} \cup \{0\}$ by  
\[ c(u) = |\{x \in V(T)|x_M=u\}|  \]
for $u \in V(T)$. Note that $c(u) =0$ if $u \not \in V(M)$ and $\sum_{u\in V(M)} c(u) = n$, 
where $n=|V(G)|$. Since $G$ is triangle-free, no two incident vertices of an edge in $M$ have 
a common neighbour. Hence, ${\rm deg}_T(u) \geq \delta$ implying that $c(u)\geq \delta$ for 
$u\in V(M)-\{v_1\}$  and ${\rm deg}_T(v_1) = \Delta$.

Now $d_T(x,x_M)\leq 3$ for every vertex $x$ of $T$. The same argument as in the proof 
of Theorem \ref{thm5-1} (see \eqref{eqn-5b}) shows that 
\begin{equation}\label{eqn5-2a}
{\rm avec}(T) \leq {\rm avec_c}(T) +3.
\end{equation} 
Let $\bar{c}$ be the weight function on the vertex set $E(T)$ of the line graph $L= L(T)$ 
defined by:
\[
\bar{c}(uv) =
\begin{cases}
c(u) +c(v) & \text{if } uv \in M, \\
0
& \text{if } uv \not \in M.
\end{cases}
\] 
Since $e_1$ is an edge incident with $v_1$ in $T$, we have $\bar{c}(e_1)\geq \Delta +\delta$ 
and $\bar{c}(e) \geq 2\delta$ for $e\in M-\{e_1\}$.  
Note that $\sum_{e\in M}\bar{c}(e) = \sum_{v\in V(T)}{c}(v)=n.$ It follows that
\begin{equation*}
\begin{aligned}
n &\geq \Delta +\delta + \sum_{x\in M-\{e_1\}} 2\delta = \Delta +\delta +2\delta (|M|-1),
\end{aligned}
\end{equation*}
and rearranging yields
\begin{equation}\label{eqn5-m}
|M| \leq \frac{n-\Delta +\delta}{2\delta}.
\end{equation}
It is easy to show that if $u, v\in V(T)$ and $e_u, e_v$ are edges of $T$ incident with $u$ and $v$, respectively, then
\[|d_L(e_u,e_v) - d_T(u,v)|\leq 1. \] 
Hence, if $u$ is an eccentric vertex of $v$, i.e., a vertex with $d_T(v,u)=e_T(v)$, then 
\[ e_T(v) = d_T(u,v)\leq d_L(e_u,e_v) +1 \leq e_L(e_v) +1,  \] 
where $e_u$ and $e_v$ are edges of $T$ incident with $u$ and $v$, respectively.
Summation over all vertices of $T$ yields that  
\begin{equation*}
\begin{aligned}
\sum_{v\in V(T)}{e_T}(v) c(v) &= \sum_{uv\in M}\big[e_T(u)c(u) + e_T(v)c(v)\big]\\
&\leq \sum_{uv \in M}\Big[\bar{c}(uv)\big(e_L(uv) +1\big)\Big]\\
&= \sum_{e \in M}e_L(e)\bar{c}(e) + \sum_{e \in M}\bar{c}(e).
\end{aligned}
\end{equation*} 
Therefore,
\[ EX_{c}(T) \leq EX_{\bar{c}}(L)+ \sum_{e \in M}\bar{c}(e),\]
and since $\sum_{v\in V(T)} c(v) = \sum_{e \in M} \overline{c}(e) = n$, 
division by $n$ yields 
\begin{equation}\label{eqn5-2b}
{\rm avec_c}(T) \leq {\rm avec}_{\bar{c}}(L) +1.
\end{equation}
If $e_1$ and $e_2$ are two matching edges with $d_T(e_1, e_2)=3$, then $d_L(e_1, e_2) \leq 4$. 
Now the weights lie solely on $M$. Thus we construct the induced subgraph $L^4[M]$ which is 
connected and has exactly $|M|$ vertices. 
Clearly, for every two edges $e_1,e_2 \in M$, $$d_L(e_1,e_2) \leq 4d_{L^4[M]}(e_1,e_2).$$ Furthermore, for every edge $e' \in E(T)$ there exists an edge $e'' \in M$ with $d_L(e',e'') \leq 3$, so that $e_L(e) \leq 4e_{L^4[M]}(e) +3$ for every $e\in M$, and so
\begin{equation}\label{eqn5-2c}
{\rm avec}_{\bar{c}}(L) \leq 4{\rm avec}_{\bar{c}}(L^4[M]) +3.
\end{equation}
We now modify the weight function $\bar{c}$ to obtain a weight function 
$\overline{c}'$ on $M$ with $\bar{c}'(e) \geq 1$ for all $e\in M$. 
Define $\bar{c}'(e) =\frac{\bar{c}(e)}{2\delta}$ for $e\in M-\{e_1\}$, and $\bar{c}'(e_1)=\frac{\bar{c}(e_1)-\Delta +\delta}{2\delta}$.
Since $\bar{c}(e_1)\geq \Delta +\delta$ and $\bar{c}(e) \geq 2\delta$ for
$e\in M-\{e_1\}$, weight function $\bar{c}'$ assigns a weight of at least $1$
to every edge of $M$. 
Let $N' = \sum_{e\in V(M)} \bar{c}'(e)$. Then $N'=\frac{n-\Delta +\delta}{2\delta}$. 
As in the proof of Theorem \ref{thm5-1} (see \eqref{eq:avec(Tc)-in-terms-of-avec(Tc')}), 
we show that  
\begin{equation} \label{L1}
avec_{\bar{c}}(L^4[M])  = 
\frac{n-\Delta+\delta}{n} 
avec_{\bar{c}'}(L^4[M]) 
+ \frac{\Delta-\delta}{n} e_{L^4[M]}(e_1).
\end{equation}

Since $L^4[M]$ has order $|M|$, and since $|M| \leq \frac{n-\Delta+\delta}{2\delta}$   
by \eqref{eqn5-m}, we have 
\begin{equation}  \label{L2} 
e_{T^4[M]}(v_1) \leq {\rm diam}(L^4[M]) \leq |M|-1 \leq \frac{n-\Delta}{2\delta} -\frac{1}{2}. 
\end{equation}
Since $ \lceil N' \rceil =   \lceil \frac{n-\Delta+\delta}{2\delta} \rceil < \frac{n-\Delta}{2\delta} + \frac{3}{2}$, we have by Proposition \ref{prop2}, 
\begin{equation}\label{L3}
{\rm avec_{\bar{c}'}}(L^4[M]) 
\leq {\rm avec}(P_{\lceil N'\rceil}) 
\leq \frac{3}{4}\lceil N'\rceil-\frac{1}{2} 
< \frac{3}{4} \frac{n - \Delta}{2\delta} + \frac{5}{8}.
\end{equation}

Substituting  \eqref{L2} and \eqref{L3} into 
\eqref{L1} yields, after simplification, 
\begin{eqnarray}
{\rm avec_{\bar{c}}}(L^4[M]) &\leq & 
\frac{n-\Delta+\delta}{n} \Big( \frac{3}{4} \frac{n - \Delta}{2\delta} + \frac{5}{8} \Big) 
+\frac{\Delta -\delta}{n} \Big( \frac{n-\Delta}{2\delta} -\frac{1}{2} \Big) \nonumber \\ 
& = & \frac{3(n-\Delta)}{8\delta} 
+ \frac{(n-\Delta)(\Delta-\delta)}{8n\delta} - \frac{9(\Delta-\delta)}{8n} +\frac{5}{8}.
\label{eqn5-3}
\end{eqnarray}
Hence combining these inequalities \eqref{eqn5-2a}, \eqref{eqn5-2b}, \eqref{eqn5-2c} and \eqref{eqn5-3} yields
\begin{eqnarray*}
{\rm avec}(T) &\leq & {\rm avec_c}(T) +3 \\
& \leq & {\rm avec}_{\bar{c}}(L) +4 \\
&\leq 4 & {\rm avec}_{\bar{c}}(L^4[M]) +7\\
&\leq & 4\big( \frac{3(n-\Delta)}{8\delta} 
+ \frac{(n-\Delta)(\Delta-\delta)}{8n\delta} - \frac{9(\Delta-\delta)}{8n} +\frac{5}{8} \Big) +7\\
&=&  \frac{n-\Delta}{2\delta} 
    \big( 3 + \frac{\Delta-\delta}{n} \big) 
    - \frac{9(\Delta-\delta)}{2n} 
    + \frac{19}{2} \\
& \leq &  \frac{3}{2} \frac{n-\Delta}{\delta} 
    \big( 1 + \frac{\Delta-\delta}{3n} \big) 
    + \frac{19}{2}. 
\end{eqnarray*}
The bound now follows since ${\rm avec} (G) \leq {\rm avec} (T)$.

\begin{figure} 
		\setlength{\unitlength}{2.3cm}
		\begin{picture}(25,2)(-3.3,-0.2)\thicklines
		\put(-3,1){\circle*{0.1}}
		\put(-3,1.5){\circle*{0.1}}
		\put(-2.5,1.5){\circle*{0.1}}
		\put(-2.5,1){\circle*{0.1}}
		\put(-2.5,2){\circle*{0.1}}

		\put(-2.5,2.4){\circle*{0.1}}
		\put(-2.5,2.8){\circle*{0.1}}
		\put(-2.5,0.0){\circle*{0.1}}
		\put(-2.5,-0.5){\circle*{0.1}}
		
		\qbezier(-3,1.5)(-2.5,2.4)(-2.5,2.4)
		\qbezier(-3,1.5)(-2.5,2.8)(-2.5,2.8) 
		\qbezier(-3,1.5)(-2.5,0.0)(-2.5,0.0)
		\qbezier(-3,1.5)(-2.5,-0.5)(-2.5,-0.5)
		
		\qbezier(-3,1.0)(-2.5,2.4)(-2.5,2.4)
		\qbezier(-3,1.0)(-2.5,2.8)(-2.5,2.8) 
		\qbezier(-3,1.0)(-2.5,0.0)(-2.5,0.0)
		\qbezier(-3,1.0)(-2.5,-0.5)(-2.5,-0.5)
		
		\qbezier(-2.0,1.5)(-2.5,2.4)(-2.5,2.4)
		\qbezier(-2.0,1.5)(-2.5,2.8)(-2.5,2.8)
		\qbezier(-2.0,1.5)(-2.5,0.0)(-2.5,0.0)
		\qbezier(-2.0,1.5)(-2.5,-0.5)(-2.5,-0.5)
		
		\qbezier(-2.0,1.0)(-2.5,2.4)(-2.5,2.4)
		\qbezier(-2.0,1.0)(-2.5,2.8)(-2.5,2.8)
		\qbezier(-2.0,1.0)(-2.5,0.0)(-2.5,0.0)
		\qbezier(-2.0,1.0)(-2.5,-0.5)(-2.5,-0.5)
		
		\put(-2.5,0.5){\circle*{0.1}}
		\put(-2.0,1.5){\circle*{0.1}}
		\put(-2.0,1.0){\circle*{0.1}}
		\put(-1.5,1.5){\circle*{0.1}}
		\put(-1.5,1.0){\circle*{0.1}}
		\put(-1.0,1.5){\circle*{0.1}}
		\put(-1.0,1.0){\circle*{0.1}}
		\put(-0.5,1.5){\circle*{0.1}}
		\put(-0.5,1.0){\circle*{0.1}}
		\put(0,1.5){\circle*{0.1}}
		\put(0,1.0){\circle*{0.1}}
		\put(0.5,1.5){\circle*{0.1}}
		\put(0.5,1.0){\circle*{0.1}}
		\put(1,1.5){\circle*{0.1}}
		\put(1,1.0){\circle*{0.1}}
		\put(1.5,1.5){\circle*{0.1}}
		\put(1.5,1.0){\circle*{0.1}}
		\put(2,1.5){\circle*{0.1}}
		\put(2,1.0){\circle*{0.1}}
		\put(2,2){\circle*{0.1}}
		\put(2,0.5){\circle*{0.1}}
		\put(2.5,1.5){\circle*{0.1}}
		\put(2.5,1.0){\circle*{0.1}}
		\qbezier(-3,1)(-2.5,1)(1,1)
		\qbezier(1,1)(1,1)(2.5,1)
		\qbezier(-3,1.5)(0,1.5)(2.5,1.5)
		\qbezier(-3,1.5)(-2.5,2)(-2.5,2)
		\qbezier(-3,1.5)(-2.5,1.0)(-2.5,1.0)
		\qbezier(-3,1.5)(-2.5,0.5)(-2.5,0.5)
		
		\qbezier(-3,1.0)(-2.5,2)(-2.5,2)
		\qbezier(-3,1.0)(-2.5,0.5)(-2.5,0.5)
		\qbezier(-3,1.0)(-2.5,1.5)(-2.5,1.5)
		
		\qbezier(-2.5,2)(-2.0,1.5)(-2.0,1.5)
		\qbezier(-2.5,1.5)(-2.0,1.0)(-2.0,1.0)
		\qbezier(-2.5,2)(-2.0,1.0)(-2.0,1.0)
		\qbezier(-2.5,1)(-2.0,1.5)(-2.0,1.5)
		\qbezier(-2.5,0.5)(-2.0,1.5)(-2.0,1.5)
		\qbezier(-2.5,0.5)(-2.0,1.0)(-2.0,1.0)
		\qbezier(-1.5,1)(-2.0,1.5)(-2.0,1.5)
		\qbezier(-2.0,1)(-1.5,1.5)(-1.5,1.5)
		\qbezier(-1.5,1)(-1.0,1.5)(-1.0,1.5)
		\qbezier(-1.5,1.5)(-1.0,1.0)(-1.0,1.0)
		\qbezier(-0.5,1.5)(-1.0,1.0)(-1.0,1.0)
		
		\qbezier(-1.0,1.5)(-0.5,1.0)(-0.5,1.0)
		\qbezier(-1.5,1.5)(-1.0,1.0)(-1.0,1.0)
		
		\qbezier(0,1.5)(-0.5,1.0)(-0.5,1.0)
		\qbezier(-0.5,1.5)(0,1)(0,1)
		
		\qbezier(0,1.5)(0.5,1.0)(0.5,1.0)
		\qbezier(0.5,1.5)(0,1)(0,1)
		
		\qbezier(0.5,1)(1.0,1.5)(1.0,1.5)
		\qbezier(0.5,1.5)(1.0,1)(1.0,1)
		
		\qbezier(1,1)(1.5,1.5)(1.5,1.5)
		\qbezier(1,1.5)(1.5,1)(1.5,1)
		
		\qbezier(1.5,1)(2,1.5)(2,1.5)
		\qbezier(2,1)(1.5,1.5)(1.5,1.5)
		
		\qbezier(2.5,1)(2,1.5)(2,1.5)
		\qbezier(2.5,1.5)(2,1)(2,1)
		
		\qbezier(2.5,1)(2,0.5)(2,0.5)
		\qbezier(2.5,1.5)(2,2)(2,2)
		
		\qbezier(2.5,1)(2,2)(2,2)
		\qbezier(2.5,1.5)(2,0.5)(2,0.5)
		
		\qbezier(1.5,1.5)(2,2)(2,2)
		\qbezier(1.5,1.0)(2,0.5)(2,0.5)
		
		\qbezier(1.5,1.0)(2,2)(2,2)
		\qbezier(1.5,1.5)(2,0.5)(2,0.5)
		\end{picture}
\caption{$G'_{\delta,\Delta,k}$ with $\Delta=10$, $\delta=4$ and $k=12$.}
\label{fig:sharpness-example-triangle-free}
\end{figure}
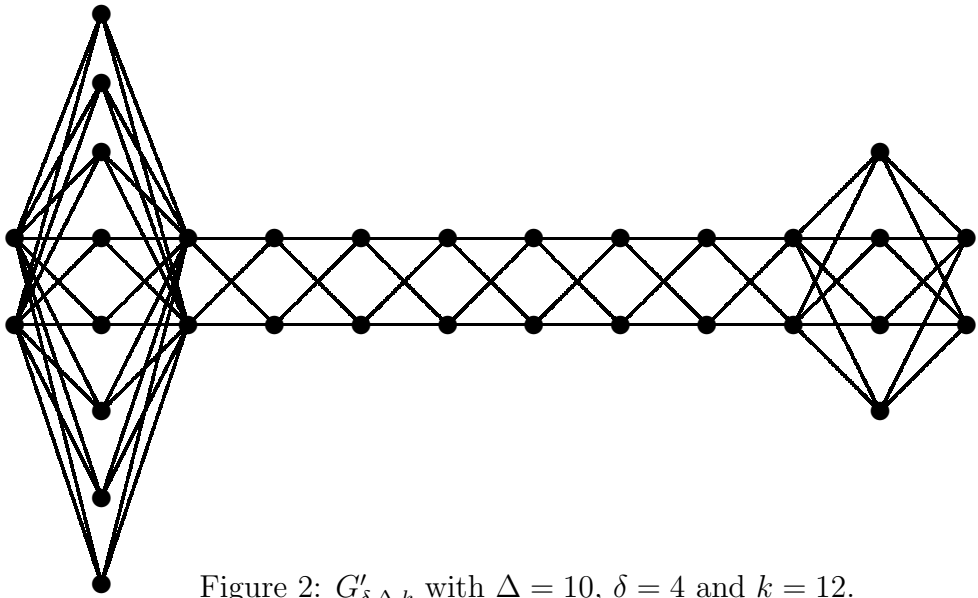

To see that the bound is sharp apart from an additive constant, let $n, \delta, \Delta$ and $k$ be positive integers with $\Delta \geq \frac{3}{2}\delta$ and $k$ a multiple of 4. Consider the graph $G'= G_{\delta, \Delta,k}'$ defined by 
\[G'=G_1+G_2+ \ldots +G_k,   \] where $G_2= (\Delta -\lceil \frac{\delta}{2}\rceil) K_1$, $G_{k-1} = \delta K_1$, and for $1\leq i\leq k$, $i \neq 2,k-1$,  
\[
G_i =
\begin{cases}
\lceil \frac{\delta}{2}\rceil K_1 & \text{if}\,\, i=0\,\, \text{or}\,\, 3\,\mod 4,  \\
\lfloor \frac{\delta}{2}\rfloor K_1 & \text{if} \,\, i=1\,\, \text{or}\,\, 2 \,\mod 4.
\end{cases}
\] 
The graph $G'_{10,4,12}$ is shown in Figure \ref{fig:sharpness-example-triangle-free}. 

By inspection, $G_{\delta, \Delta,k}'$ has minimum degree $\delta$, maximum degree $\Delta$, 
and it is triangle-free. Moreover, $G_{\delta, \Delta,k}'$ has exactly 
$\Delta + \delta (\frac{k}{2}) -\lceil \frac{\delta}{2}\rceil$ vertices. A simple 
calculation, whose details we omit, yields that 
\begin{equation*}
EX(G'_{\delta,\Delta,k}) 
= \delta\Big( \frac{3k^2}{8}-\frac{k}{4}\Big)
     +(\Delta -\Big\lceil \frac{\delta}{2}\Big\rceil)(k-2).
\end{equation*}
Let $n= n(G'_{\delta,\Delta,k}) 
=\Delta + \delta (\frac{k}{2}) -\lceil \frac{\delta}{2}\rceil$. Then 
$k=\frac{2}{\delta}(n-\Delta+\lceil \frac{\delta}{2} \rceil) 
      = \frac{2}{\delta}(n-\Delta) + O(1)$. 
Hence 
$\delta\Big( \frac{3k^2}{8}-\frac{k}{4}\Big) 
 =  \frac{3}{8}\delta k^2 +O(n)
 = \frac{3}{2\delta} (n-\Delta)^2 + O(n)$
and
$ (\Delta -\Big\lceil \frac{\delta}{2}\Big\rceil)(k-2)
  = \Delta k +O(n) 
  = \frac{2}{\delta}(n-\Delta)\Delta + O(n)$. 
Substituting this into the above expression for  $EX(G'_{\delta,\Delta,k})$
and dividing by $n$ we get
\begin{eqnarray*}
{\rm avec}(G_{\delta,\Delta,k}) 
   &=&  \frac{1}{n} 
    \Big( \frac{3}{2\delta} (n-\Delta)^2   +  \frac{2}{\delta}(n-\Delta)\Delta + O(n) \Big) \\
    &=& \frac{3}{2} \frac{n-\Delta}{\delta} \big(1 + \frac{\Delta}{3n}\big) + O(1). 
\end{eqnarray*}
Hence the bund in Theorem \ref{thm5-2} and ${\rm avec}(G_{\delta,\Delta,k})$ is
bounded by a constant.   
\qed

\subsection{A bound on the average eccentricity of a $C_4$-free graph}

\noindent
In this section we show that the bound in Theorem 1 can be improved for graphs not 
containing $C_4$ as a (not necessarily induced) subgraph.

\begin{theorem}\label{theo:C4-free}
Let $G$ be a $C_4$-free graph of order $n$, minimum degree $\delta$ and maximum degree 
$\Delta$. Then 
\[{\rm avec}(G)\leq
\frac{15}{4} \frac{n-\varepsilon_{\Delta} + \varepsilon_{\delta}}{\varepsilon_{\delta}} 
\Big[ 1 + \frac{\varepsilon_{\Delta} - \varepsilon_{\delta}}{3n} \Big]  + \frac{37}{4},
	\]
where 
$\varepsilon_{\Delta}:= \Delta \delta -2\Big\lfloor\frac{\Delta}{2}\Big\rfloor +1$
and 
$\varepsilon_{\delta} := \delta^2 -2\lfloor \frac{\delta}{2} \rfloor +1$. 		
\end{theorem}

\noindent  
{\bf Proof.} Let $v_1$ be a vertex in $G$ with ${\rm deg}_G(v_1)=\Delta$. We find a $4$-packing $A$ of $G$ using the following method. Let $A=\{v_1\}$. If there exists a vertex $a_1$ with $d_G(a_1, A)=5$, add $a_1$ to $A$. Add vertices $a_i$ satisfying $d_G(a_i,A)=5$ to $A$ until each of the vertices not in $A$ is within distance $4$ of $A$. For each $a\in A$, let $T_{a}$ be a subtree of $G$ with vertex set $N_G^2[a]$ which is distance preserving from $a$. Let $T_1= \bigcup_{a\in A} T_{a}$. Then $T_1$ is a subforest of $G$. By the way $A$ is constructed, there exists $|A|-1$ edges in $G$, each joining two components of $T_1$, whose addition to $T_1$ yields a subtree $T_2$ of $G$.

Let $T$ be a spanning tree of $G$ obtained from $T_2$ and which satisfies $d_T(x,A) =d_G(x,A)$ for each $x\in V(G)$. 
Since ${\rm deg}_T(v_1) = {\rm deg}_G(v_1)$, tree $T$ has the same maximum degree as $G$. 
For every vertex $x$ of $T$ let $x_A$ be a vertex in $A$ closest to $x$. 
We define a weight function $c: V(T)\rightarrow \mathbb{N}\cup \{0\}$ for each $u\in V(T)$ as 
\[c(u) =|\{x\in V(T)|x_A=u \} |.\]
Since $G$ is $C_4$-free, no two neighbours of $v_1$ have a common neighbour apart from $v_1$. It follows that for $\Delta$ even, we have
$|N_G^2[v_1]| \geq 1 + \Delta + \Delta(\delta -2)=\Delta \delta -\Delta +1$.  
For $\Delta$ odd, the handshake lemma yields that at least one of the neighbours 
of $v_1$ is not adjacent to any other neighbour of $v_1$, and so 
$|N_G^2[v_1]| \geq 1 + \Delta + (\Delta-1)(\delta -2) +(\delta -1) =\Delta \delta -\Delta +2$. 
Combining these two bounds on $|N_G^2[v_1]|$ yields that 
$|N_G^2[v_1]| \geq \Delta \delta -2\Big\lfloor\frac{\Delta}{2}\Big\rfloor +1
     = \varepsilon_{\Delta}$. 
The same reasoning as for $v_1$ shows that for vertices $a_i \in A-\{v_1\}$, 
$|N_G^2[a_i]| \geq \delta^2 -2\lfloor \frac{\delta}{2} \rfloor +1 =\varepsilon_{\delta}$. 
This implies that 
\begin{equation}  \label{eq:lower-bound-on-c} 
c(v_1) \geq \varepsilon_{\Delta}, \quad 
\textrm{and $c(a_i) \geq \varepsilon_{\delta}$ for $a_i\in A-\{v_1\}$},  
\end{equation}
while $c(u)=0$ whenever $u\not \in A$. 
Now for every vertex $x \in V(T)$, $d(x,x_A) \leq 4$. Thus we have  
$|e_T(x) - e_T(x_A)| \leq 4$, and an argument similar to \eqref{eqn-5b} shows that 
\begin{equation}\label{eqn5.1.3c}
{\rm avec}(T)\leq {\rm avec_c}(T) +4.
\end{equation}
By the way $A$ was constructed, $T^5[A]$ is connected and so for $a_i,a_j\in A$, $d(a_i,a_j)\leq 5d_{T^5[A]}(a_i,a_j)$. Since every vertex of $T$ is within distance 4 of $A$, we have
\begin{equation}\label{eqn5.1.3d}
{\rm avec_c}(T) \leq 5{\rm avec_c}(T^5[A]) +4.
\end{equation}
We modify $c$ to obtain a weight function $c''$ for which every vertex of $A$ has weight
at least $1$. 
Define $c''(a_i) = \frac{c(a_i)}{\varepsilon_{\delta}}$ 
for all $a_i \in A-\{v_1\}$ and 
$c''(v_1)= 
\frac{c(v_1)-\varepsilon_{\Delta} + \varepsilon_{\delta}}{\varepsilon_{\delta}}$. 
By \eqref{eq:lower-bound-on-c} 
we have  $c''(a) \geq 1$ for all $a\in A$. 
 
Let $N'' = \sum_{u\in A}c''(u)$. Then 
\[ N''= \frac{n-\varepsilon_{\Delta} + \varepsilon_{\delta}}{\varepsilon_{\delta}}.\] 
As in the proof of Theorem 1 (see \eqref{eq:avec(Tc)-in-terms-of-avec(Tc')}),
we express ${\rm avec_{c}}(T^5[A])$ in terms of ${\rm avec_{c''}}(T^5[A])$. Clearly,  
\begin{eqnarray*}
{\rm avec_{c''}}(T^5[A]) 
&=&  \frac{EX_{c''}(T^5[A])}{N''} \\
&=& \frac{ \frac{1}{\varepsilon_{\delta}}  EX_{c}(T^5[A])
	  - \frac{\varepsilon_{\Delta} - \varepsilon_{\delta}}{\varepsilon_{\delta}}  
	   e_{T^5[A]}(v_1)}{(n-\varepsilon_{\Delta}+ \varepsilon_{\delta}) / \varepsilon_{\delta}}  \\ 
&=& \frac{n}{n-\varepsilon_{\Delta} +\varepsilon_{\delta}} {\rm avec_c}(T^5[A]) 
  - \frac{\varepsilon_{\Delta} - \varepsilon_{\delta}}{n-\varepsilon_{\Delta} +\varepsilon_{\delta}} e_{T^5[A]}(v_1).  
\end{eqnarray*}
By rearranging, we obtain 
\begin{equation}\label{t5}
{\rm avec_c}T^5[A] 
 = \frac{n-\varepsilon_{\Delta} + \varepsilon_{\delta}}{n} {\rm avec_{c''}}(T^5[A]) 
   + \frac{\varepsilon_{\Delta} - \varepsilon_{\delta}}{n} e_{T^5[A]}(v_1).
\end{equation} 
We bound the terms on the right hand side of \eqref{t5} separately. 
Since $T^5[A]$ has order $|A|$, we have $e_{T^5[A]}(v_1) \leq |A|-1$. 
Also $|A| = \sum_{a\in A} 1 \leq  \sum_{a\in A} c''(a) = N''$, and so  
 \begin{equation}\label{t6}
 e_{T^5[A]}(v_1) \leq |A|-1  \leq N'' - 1
   = \frac{n-\varepsilon_{\Delta}+\varepsilon_{\delta}}{\varepsilon_{\delta}} - 1. 
 \end{equation}
Since $\lceil N'' \rceil < N''+1 
= \frac{n - \varepsilon_{\Delta} +\varepsilon_{\delta}}{\varepsilon_{\delta}}  
  + 1$ 
it follows by Proposition \ref{prop2} that 
\begin{equation}\label{eqn5.1.3g}
{\rm avec_{c''}}(T^5[A]) 
\leq {\rm avec}(P_{\lceil N''\rceil}) 
\leq \frac{3\lceil N''\rceil}{4}-\frac{1}{2}
<  \frac{3}{4} \frac{n-\varepsilon_{\Delta} + \varepsilon_{\delta}}{\varepsilon_{\delta}} 
  +\frac{1}{4}.
\end{equation}
Thus substituting \eqref{eqn5.1.3g} and \eqref{t6} in \eqref{t5} yields, after simplification,
\begin{eqnarray*}
{\rm avec_c}T^5[A] 
 & < &\frac{n-\varepsilon_{\Delta} + \varepsilon_{\delta}}{n} 
  \Big[ \frac{3}{4} \frac{n-\varepsilon_{\Delta} + \varepsilon_{\delta}}{\varepsilon_{\delta}} 
  +\frac{1}{4} \Big] 
   + \frac{\varepsilon_{\Delta} - \varepsilon_{\delta}}{n} 
   \big( \frac{n-\varepsilon_{\Delta}+\varepsilon_{\delta}}{\varepsilon_{\delta}}
      - 1 \Big) \nonumber \\
& = & \frac{3}{4} \frac{n - \varepsilon_{\Delta} + \varepsilon_{\delta}}{ \varepsilon_{\delta}} 
      \Big[ 1 + \frac{ \varepsilon_{\Delta} -\varepsilon_{\delta}}{3n} \Big] + \frac{1}{4}  
        - \frac{5}{4} \frac{\varepsilon_{\Delta}-\varepsilon_{\delta}}{n}  \nonumber  \\
& \leq & \frac{3}{4} \frac{n - \varepsilon_{\Delta} + \varepsilon_{\delta}}{ \varepsilon_{\delta}} 
      \Big[ 1 + \frac{ \varepsilon_{\Delta}-\varepsilon_{\delta}}{3n} \Big] + \frac{1}{4}.    \label{t7}        
\end{eqnarray*} 
Combining inequalities \eqref{eqn5.1.3c}, \eqref{eqn5.1.3d} and \eqref{t7}, we obtain with further simplification,
\begin{eqnarray*}
{\rm avec}(T) & \leq &{\rm avec_c}(T) +4 \\
&\leq &  5 {\rm avec_{c}}(T^5[A]) +8 \\
   & \leq & \frac{15}{4}  \frac{n- \varepsilon_{\Delta} + \varepsilon_{\delta}}{\varepsilon_{\delta}} 
         \big[1 + \frac{\varepsilon_{\Delta} - \varepsilon_{\delta}}{3n} \big]  
          + \frac{37}{4},    
\end{eqnarray*}
as desired.
\qed   \\

The following theorem shows that the bound in Theorem \ref{theo:C4-free} is not far
from being sharp if $\delta +1$ is a prime power.  In our construction below, the 
maximum degree and order can be chosen almost arbitrarily, they only have to satisfy 
certain conditions regarding their values modulo $\delta+1$ or $\delta+2$.   
Our construction is a modification of a graph first constructed  by 
Erd\"{o}s et al.\ \cite{ErdPacPolTuz1989}.

\begin{theorem}
Let $\delta \geq 3$ be an integer such that $\delta +1$ is a prime power. 
Then for $n, \Delta \in \mathbb{N}$ with
$2\delta-3 \leq \Delta <n$ and $n \equiv 0 \pmod{(\delta+1)(\delta+2)}$ 
and $\Delta \equiv \delta +1 \pmod{\delta+2}$ there exists a $C_4$-free graph 
$G$ with minimum degree $\delta$, maximum degree $\Delta$, order
$n$ whose average eccentricity satisfies 
\[{\rm avec}(G)\geq 
 \frac{3}{4} \frac{n - \varepsilon_{\Delta} - \Delta}{\varepsilon_{\delta}'} 
    \big( 1 + \frac{\Delta (\delta+1)}{3n} \big) + O(1), \]
where $\varepsilon_{\Delta}$ is as in Theorem \ref{theo:C4-free},
and $\varepsilon_{\delta}' =(\delta+1)(\delta+2)$.       
\end{theorem}

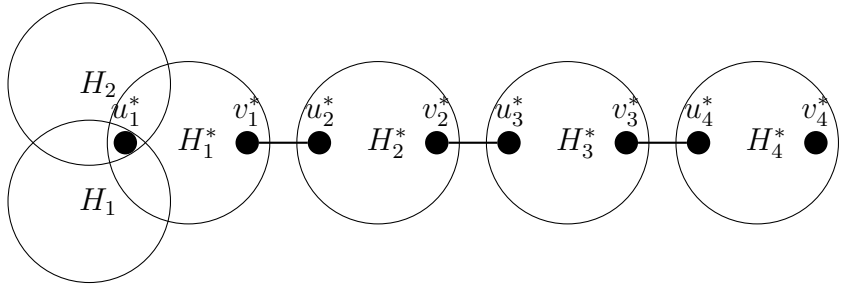
\begin{figure}
  \begin{center}
\begin{tikzpicture}
  [scale=0.6,inner sep=1mm, 
   vertex/.style={circle,thick,draw}, 
   thickedge/.style={line width=2pt}] 

    \node[vertex] (a1) at (0.8,1.3) [fill=black] {};
    \node[vertex] (a2) at (3.5,1.3) [fill=black] {};    
    \node[vertex] (a3) at (5.1,1.3) [fill=black] {};
    \node[vertex] (a4) at (7.7,1.3) [fill=black] {};
    \node[vertex] (a5) at (9.3,1.3) [fill=black] {};    
    \node[vertex] (a6) at (11.9,1.3) [fill=black] {};
    \node[vertex] (a7) at (13.5,1.3) [fill=black] {};   
    \node[vertex] (a8) at (16.1,1.3) [fill=black] {};     
        
\draw[black] (0,0) circle [radius=1.8];
\draw[black] (0,2.6) circle [radius=1.8];
\draw[black] (2.2,1.3) circle [radius=1.8];
\draw[black] (6.4,1.3) circle [radius=1.8];
\draw[black] (10.6,1.3) circle [radius=1.8];
\draw[black] (14.8,1.3) circle [radius=1.8];

\draw[thick] (a2)--(a3)  (a4)--(a5)   (a6)--(a7);

\node[right] at (0.3,2) {$u_1^*$};
\node[right] at (3,2) {$v_1^*$};
\node[right] at (4.6,2) {$u_2^*$};
\node[right] at (7.2,2) {$v_2^*$};
\node[right] at (8.8,2) {$u_3^*$};
\node[right] at (11.4,2) {$v_3^*$};
\node[right] at (13,2) {$u_4^*$};
\node[right] at (15.6,2) {$v_4^*$};

\node at (2.4,1.3) {$H_1^*$};
\node at (6.6,1.3) {$H_2^*$};
\node at (10.8,1.3) {$H_3^*$};
\node at (15,1.3) {$H_4^*$};
\node at (0.2,0.0) {$H_1$};
\node at (0.2,2.6) {$H_2$};
\end{tikzpicture}
\end{center}
\caption{$H_{\delta,\Delta,k,m}$ for $k=4$ and $m=2$.}
\label{fig:C4-free-sharpness-example}
\end{figure}

\noindent
{\bf Proof:} Our construction is a modification of a construction given 
in \cite{ErdPacPolTuz1989}. 
Let $q=\delta+1$, so $q$ is a prime power. Denote be $GF(q)$ the field of 
order $q$ and by $GF(q)^3$ the 
$3$-dimensional vector space over $GF(q)$ of all triples of elements of $GF(q)$. 
Let $H$ be the graph whose vertices are the $1$-dimensional subspaces of $GF(q)^3$. 
Two vertices are adjacent if, as subspaces of $GF(q)^3$, they are orthogonal. 
It is easy to verify that $H$ is $C_4$-free, has $q^2+q+1$ vertices and that 
${\rm diam}(H)=2$. Each vertex of $H$ 
has degree either $q+1$ (if the corresponding subspace is not self-orthogonal) or 
$q$ (if the corresponding subspace is self-orthogonal). 

Let $H^*$ be the graph obtained from $H$ by choosing a vertex $z$ that corresponds to
a self-orthogonal subspace, two neighbours $u$ and $v$ of $z$, and deleting vertex
$z$ as well as all edges joining a neighbour of $u$ to a neighbour of $v$. It is 
easy to show that $H^*$ has order $q^2+q$, that $\delta(H^*) \geq q-1$, and that
$d_{H^*}(u,v)= {\rm diam}(H^*)=4$.

For $i=1,2,\ldots,k$ let $H^*_i$ be a copy of $H^*$. We denote the vertices of $H^*_i$ 
corresponding to the vertices $u$ and $v$ of $H^*$ by $u_i^*$ and $v_i^*$, respectively. 
For $j=1,2\ldots,\ell$ let $H_i$ be a copy of $H$, and let $w_i$ be a vertex of 
$H_i$ of degree $q+1$.  Let $H_{\delta,\Delta,k,\ell}$ be the graph obtained
from the disjoint union 
$H^*_1 \cup H^*_2 \cup \ldots,\cup H^*_k \cup H_1\cup H_2\cup \ldots \cup H_{\ell}$
by adding the edges $v_1^*u_2^*, v_2^*u_3^*,\ldots,v_{k-1}^*u_k^*$ and by identifying the 
vertices $w_1, w_2,\ldots,w_{\ell}$ and $u_1^*$ to a new vertex $y$. 
Then $H_{\delta,\Delta,k,\ell}$ has order $n=(m+k)(q^2+q)$, minimum degree
$\delta=q-1$, and maximum degree $\Delta=(m+1)(q+1)-1$, which is attained by vertex $y$.  

We now bound the average eccentricity of this graph from below. 
For $i=0,1,\ldots,\frac{k}{2}-1$, each of the $2(q^2+q)$ vertices in 
$V(H^*_{k/2-i}) \cup V(H^*_{k/2+1+i})$ has eccentricity at least $5(i+k/2)$. 
The remaining $m(q^2+q)$ vertices, which are in $(V(H_1) \cup \ldots V(H_m))-\{u_1^*\}$, 
have eccentricity at least $5k$. Hence, 
\begin{eqnarray*} 
EX(H_{\delta,\Delta,k,m}) 
  &=& \big( \sum_{i=0}^{k/2-1} 2(q^2+q) 5 (\frac{k}{2}+i) \big) + m(q2+q) 5k  \\
  &=& 5k(q^2+q) \big( \frac{3}{4}k + m - \frac{1}{2} \big). 
\end{eqnarray*}  
We now express the two factors of the above term for $EX(H_{\delta,\Delta,k,m})$ 
separately in terms of $n$, $\Delta$ and $\delta$. 
Let $\varepsilon_{\Delta}$ and $\varepsilon_{\delta}$ be as in Theorem \ref{theo:C4-free},
and define $\varepsilon_{\delta}' =(\delta+1)(\delta+2)$. 

Since $q=\delta+1$ and $q$ is constant, we have 
$n-\varepsilon_{\Delta} + \varepsilon_{\delta} -2\Delta 
 = n - \Delta \delta +O(1)
 = (m+k) (q^2+q) - (m+1) - [(m+1)(q+1)-1](q-2) - 2(m+1)(q+1) - 2 +O(1) 
 =k(q^2+q) +O(1)$. Hence the first factor is 
 $k(q^2+q) = n-\varepsilon_{\Delta}+\varepsilon_{\delta} - 2\Delta +O(1)$. 
 Now consider the second factor. 
$\frac{3}{4}k + m -\frac{1}{2} = \frac{3}{4}(k+m) + \frac{1}{4}m +O(1) 
   = \frac{3}{4} \frac{n}{q^2+q} + \frac{1}{4} \frac{\Delta}{q+1} +O(1)
   = \frac{3}{4} \big( \frac{n}{\delta+1)(\delta+2)} + \frac{\Delta(\delta+1)}{3(\delta+1)(\delta+2)} \big) +O(1)$. Combining these terms we obtain 
\[ EX(H_{\delta, \Delta,k,m}) 
  = \frac{3}{4} \frac{n - \varepsilon_{\Delta} - \Delta}{\varepsilon_{\delta}'} 
    \big( n + \frac{\Delta (\delta+1)}{3} \big) + O(n). \]
Dividing by $n$ now yields the statement of the theorem. 
\qed

\end{document}